\documentclass[11pt,a4paper]{article}
\usepackage[a4paper, total={14.5cm,25cm}]{geometry}

\usepackage[round]{natbib}
\usepackage{sources/latex_packages/myLatexStyle/myColours}
\usepackage{sources/latex_packages/myLatexStyle/myPackagesImagePlot}
\usepackage{sources/latex_packages/myLatexStyle/myPackages}
\usepackage{sources/latex_packages/myLatexStyle/myMathOperators}
\usepackage{sources/latex_packages/myLatexStyle/myCommands}
\usepackage{sources/latex_packages/myLatexStyle/myLatexStyleArticleNeutral/myStyle}
\usepackage{sources/latex_packages/myLatexStyle/myCrefName}

\newcommand{\pathSources}{sources/}

\newcommand{\pathLatexFiles}{latex_files/}

\usepackage{\pathLatexFiles/notations}

\newcommand{\orcid}[1]{\href{https://orcid.org/#1}{\includegraphics[width=.4cm]{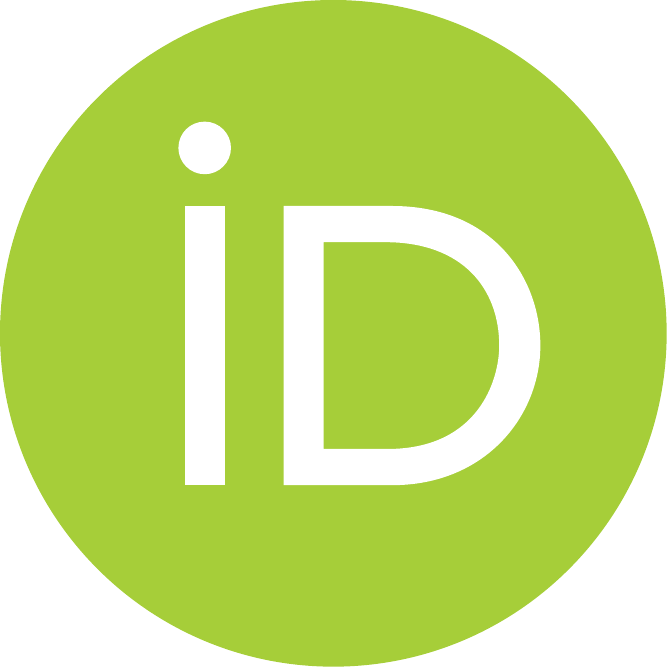}}}

\providecommand{\keywords}[1]
{
  \noindent
  \small
  \textbf{Keywords ---} #1
  \normalsize
}

\date{\monthyeardate\today}

\title{
  A data-driven uncertainty modelling and reduction approach for energy optimisation problems
}

\author{
  Julien Vaes
  \thanks{
    Department of Chemical Engineering, The Sargent Centre for Process Systems Engineering, University College London, Torrington Place, London WC1E 7JE, UK
    (\href{mailto:j.vaes@ucl.ac.uk}{j.vaes@ucl.ac.uk}).
  }
  \hspace{2mm}\orcid{0000-0003-4562-5373}%
  \and 
  Vassilis M. Charitopoulos
  \thanks{
    Department of Chemical Engineering, The Sargent Centre for Process Systems Engineering, University College London, Torrington Place, London WC1E 7JE, UK
    (\href{mailto:v.charitopoulos@ucl.ac.uk}{v.charitopoulos@ucl.ac.uk}).
  }
  \hspace{2mm}\orcid{0000-0001-9051-917X}%
}

\begin{document}

\maketitle

\begin{abstract}
  Taking uncertainty into account is crucial when making strategic decisions.
To guard against the risk of adverse scenarios,
traditional optimisation techniques incorporate uncertainty on the basis of prior knowledge on its distribution.
In this paper,
we show how,
based on a limited amount of historical data,
we can generate from a low-dimensional space the underlying structure of uncertainty that could then be used in such optimisation frameworks.
To this end,
we first exploit the correlation between the sources of uncertainty through a principal component analysis to reduce dimensionality.
Next,
we perform clustering to reveal the typical uncertainty patterns,
and finally we generate polyhedral uncertainty sets based on a kernel density estimation (KDE) of marginal probability functions.

\end{abstract}




\section{Introduction}%
\label{sec:introduction}

Optimisation under uncertainty aims at finding optimal solutions  subject to uncertain outcomes.
In most optimisation methods,
including chance-constrained optimisation, robust optimisation and stochastic optimisation,
it is assumed that an uncertainty set or a probability distribution function on these uncertain outcomes is known beforehand.
Robust optimisation is for instance a suitable method in situations where insufficient historical data is available to accurately estimate distributions \citep{Roald2023a}.
In this method,
the solution hedges against a range of possible realisations called \emph{uncertainty set} with emphasis on the worst-case to guarantee feasibility of the policy derived.
The definition of such uncertainty set is crucial:
the larger it is,
the less informative it becomes and the more conservative is the solution,
which could render it impractical \citep{Ben-Tal2009}.

We mention here some methods and their limitation that have been recently proposed in the literature to 
generate uncertainty sets for robust optimisation.
First,
the method of \citet{Bertsimas2018} employs statistical hypothesis tests to derive uncertainty sets' sizes.
Nonetheless,
in cases of data scarcity,
these tests become not significant and render the method impractical.
In \citet{Shang2017,Ning2018},
a single polyhedral uncertainty set (PUS) is obtained using respectively
kernel learning (support vector clustering) or
a principal component analysis (PCA) coupled with kernel smoothing.
This implicitly assumes that the data is formed as a continuous cloud of points and prevents the use of the method when uncertainty is formed as a collection of independent states.
To circumvent this limitation,
\citet{Ning2017,Ning2019a} define the uncertainty set as the union of PUS,
which are generated using a Dirichlet process mixture model.
In this method,
the latent uncertainty used to generate the PUS is of the same dimension of the original uncertainty,
which can be prohibitive when uncertainty belongs to a high dimensional space.
For instance,
a robust optimisation problem can be computationally intractable for highly-dimensional uncertainty sets \citet{Ben-Tal2009}.
There is consequently the need for a method that generates polyhedral uncertainty sets in the case of data sparsity where the uncertainty is highly dimensional and consists of independent sets of realisations;
this is the objective of the present paper.

The rest of this paper is structured as follows: 
In \cref{sec:method_for_deriving_polyhedral_uncertainty_and_ambiguity_sets},
we present the method which derives polyhedral uncertainty sets from a low-dimensional space based on the historical data. 
We then show in \cref{sec:numerical_example} how it can be applied to the energy data in the United Kingdom. In \cref{sec:results} discussion of the key findings is presented and finally in 
\cref{sec:conclusion} conclusions on the main contributions of this work are drawn.

\section{Method for deriving polyhedral uncertainty sets}%
\label{sec:method_for_deriving_polyhedral_uncertainty_and_ambiguity_sets}

To motivate the importance of building polyhedral uncertainty sets,
we consider the following scenario-based linear adaptive robust optimisation (ARO) problem,
which is recurrently used to formulate energy problems (\eg unit commitment):
\begin{mini!}|l|[2]
  {\xb}
  {\cb^{\T} \xb + \sum_{k \in \mset{K}} p_{k}
    \pa{
    \max_{\vb_{k} \in \mset{V}_{k}} \min_{\yb_{k}} \bb_{k}^{\T}\yb_{k}, \label{eq:RO_objective}}
  }
  {\label{op:RO_problem}}
  {}
  \addConstraint{}{\Ab \xb \leq \db, \label{eq:first_stage_constraint}}
  \addConstraint{\forall k \in \mset{K}:}{\hb_{k} - \Tb_{k}\xb - \Mb_{k} \vb_{k} \leq \Wb_{k} \yb_{k}, \label{eq:second_stage_constraint}}
\end{mini!}
where $\xb$ and $\yb_{k}$ with $k \in \mset{K}$ are respectively the first and second stage variables,
$\mset{K}$ is the set of distinct operation conditions/scenarios,
$\mset{V}_{k}$ is the uncertainty set associated to the uncertain parameters in the operation conditions $k \in \mset{K}$,
and $p_{k}$ is the weight (or probability of occurrence) associated to scenario~$k$.
Such problem can be solved either in a monolithic fashion or by using variants of Benders' decomposition.
If each set $\mset{V}_{k}$ is a polyhedral uncertainty set (PUS),
then each set can be expressed as linear constraints,
which in turn implies that the resolution of the ARO problem is tremendously facilitated by exploiting the dual formulation associated to a linear program \citep{Birge2011}.   
The rest of the section shows how to derive these polyhedral uncertainty sets~$\mset{V}_{k}$ based on limited historical data and how to generate them from lower dimensional spaces.
We should note that the method is here illustrated for energy uncertainty data but is not restricted to such application.

Let $n \in \N$ be the number of uncertain attributes (\eg electricity demand, temperature).
Let $\vb^{\spa{i}} \in \R^{24 n}$,
$i \in \rangeo{m}$,
be $m$ historical data points,
where $\vb^{\spa{i}}$ corresponds to the vector containing the standardised daily profiles (hourly values) for all the attributes on a given day $i \in \rangeo{m}$.
We assume that the data is standardised in the sense that the sample mean and standard deviation of each component over all~$m$ data points is respectively 0 and 1.
This allows for attributing an even weight to any uncertainty component.
Given a data point $\vb^{\spa{i}} \in \R^{24 n}$,
let $\Vb^{\spa{i}} \in \R^{n \times 24}$ be the matrix with the daily profiles of each attribute on its rows.

As an initial step,
based a correlation analysis,
we define a set~$\mset{A}$ of disjoint groups of attributes such that the data of the attributes inside a group are strongly correlated while being weakly correlated with any attribute of another group.
The rest of the method works as follows:
we assume that the groups of attributes are independent and separately generate a polyhedral uncertainty set for each one of them.
The (weak) correlation between the groups of attributes is then partly reinstated by defining scenarios,
which are defined as the combination of specific operation conditions for each group of attributes.

For any group of attributes $a \in \mset{A}$,
we denote by $\vb_{a}^{\spa{i}}$ the data points related to the $n_{a}$ attributes associated to~$a$ for day $i \in \rangeo{m}$.
A data point in group $a$ is thus an element of $\R^{24 n_{a}}$.
The following steps allow for deriving a PUS based on a lower dimensional space for each group $a \in \mset{A}$:
First, 
we perform a principal component analysis (PCA) and express the data in the PCA basis: 
this allows for finding the directions along which the data has greatest variance.
Mathematically,
there exists an orthogonal matrix $\Pb_{a} \in \R^{24 n_{a} \times 24 n_{a}}$ such that,
any data point $\vb_{a}^{\spa{i}}$,
$i \in \rangeo{m}$,
expressed in the original basis can be expressed in the PCA basis using the following linear transformation:
\begin{equation}
  \label{eq:pca_linear_transformation}
  \wb_{a}^{\spa{i}}
  =
  \Pb_{a} \vb_{a}^{\spa{i}}.
\end{equation}
As $\Pb_{a}$ is an orthogonal matrix,
the linear transformation from the PCA basis to the original basis is given by $ \vb_{a}^{\spa{i}} = \Pb_{a}^{\T} \wb_{a}^{\spa{i}} $.
Given the data expressed in the PCA basis,
we perform dimensionality reduction by retaining only the first~$r_{a}$ components that explain the most the variability inside the data.
As the data is strongly correlated inside any given group $a \in \mset{A}$,
we expect to have $r_{a} \ll 24 n_{a}$.
Let $\bar{\wb}_{a}^{\spa{i}} \in \R^{r_{a}}$ denote the data points in the truncated PCA basis,
where only the first $r_{a}$ components are kept.

To define the typical operational conditions associated to a group $a \in \mset{A}$ of attributes,
we perform clustering  (\eg K-means) on the data points expressed in the truncated PCA basis.
Clustering in the PCA basis rather than in the original basis allows to give relatively more importance to the directions explaining the most the variance and obtain more meaningful clusters \citep{Ding2004}.
We denote by $\mset{K}_{a}$ the set of clusters related to $a \in \mset{A}$,
call each cluster $k_{a} \in \mset{K}_{a}$ a \emph{group-scenario}
and define its probability of occurrence~$p_{k_{a}}$ as the proportion of data points attributed to this cluster.

We now follow a similar approach as \citet{Ning2018} to construct a polyhedral uncertainty set for each group-scenario $k_{a} \in \mset{K}_{a}$.
To this end,
we estimate the marginal probability density function~$\hat{f}_{a,k_{a},r}$,
$r \in \rangeo{r_{a}}$,
along each truncated principal component direction (with a kernel smoothing method such as KDE for instance \citep{Chen2017}).
We denote by~$\hat{F}_{a,k_{a},r}$ the associated cdf,
and by~$\hat{F}^{-1}_{a,k_{a},r}$ the associated quantile function.
We now define the polyhedral uncertainty set (PUS) related to the group of attributes $a \in \mset{A}$ and the cluster $k_{a} \in \mset{K}_{a}$,
which we call a \emph{group-scenario PUS},
as follows:
\begin{equation}
  \label{eq:construction_polyhedral_uncertainty_set_group}
  \bar{\mset{W}}_{a,k_{a}}^{pol}
  \mdef
  \setc{
    \bar{\wb} \in \R^{r_{a}}
    }{
    \begin{array}{l}
      \0 \leq \zb^{-}, \zb^{+} \leq \1, \\
      \1^{\T} \spa{\zb^{-} + \zb^{+}} \leq \Phi_{a,k_{a}} \\
      \forall i, j \in \rangeo{s_{a,k_{a}}}: z_{i}^{-} + z_{i}^{+} + z_{j}^{-} + z_{j}^{+} \leq \Psi_{a,k_{a}} \\
      \xib_{a,k_{a}}^{lb} = \pac{ \hat{F}_{a,1,k_{a}}^{-1} \spa{\alpha_{a,k_{a}}}, \dots, \hat{F}_{a,r_{a},k_{a}}^{-1} \spa{\alpha_{a,k_{a}}} }^{\T}  \\
      \xib_{a,k_{a}}^{ub} = \pac{ \hat{F}_{a,k_{a},1}^{-1} \spa{1-\alpha_{a,k_{a}}}, \dots, \hat{F}_{a,k_{a},r_{a}}^{-1} \spa{1-\alpha_{a,k_{a}}} }^{\T} \\
      \lambdab = \frac{1}{2} \pa{\zb^{+} - \zb^{-} + \1} \\
      \bar{\wb}
      = \xib_{a,k_{a}}^{lb} \circ \pa{1 - \lambdab} + \xib_{a,k_{a}}^{ub} \circ \lambdab
    \end{array}
  },
\end{equation}
where $\0$ and $\1$ are vectors of respectively zeros and ones of size $r_{a}$,
where $ \circ $ defines the Hadamard product,
and where vectors inequalities must be understood componentwise.
This set is similar to the budgeted uncertainty set proposed by \citep{Bertsimas2004a},
which enables to control the degree of conservatism :
the larger the volume of the PUS,
the more conservative is the resulting robust Optimisation Problem~\eqref{op:RO_problem} \citep{Bertsimas2004a}.
The uncertainty set is parametrised by $\alpha_{a,k_{a}}$, $\Phi_{a,k_{a}}$ and $\Psi_{a,k_{a}}$.
First,
$\alpha_{a,k_{a}}$
is used to exclude both tails of the marginal pdf along each principal component axis.
Then,
$\Phi_{a,k_{a}}$ limits the cumulative dispersion from the nominal value along all retained PCA axes.
Finally,
the parameter $\Psi_{a,k_{a}}$ additionally limits the pairwise dispersion along the first $s_{a,k_{a}}$ PCA components.
This parameter is important to further exclude unlikely data points from the uncertainty that would otherwise not be cut with the general budget constraint parametrised by $\Phi_{a,k_{a}}$.

We have so far generated a PUS for each group-scenario $k_{a} \in \mset{K}_{a}$ for each $a \in \mset{A}$.
Let $\mset{K}^{\times} \mdef \bigtimes_{a \in \mset{A}} \mset{K}_{a}$ denote the Cartesian product of the sets of group-scenarios.
The probability $\hat{p}_{k}$ associated to $k \in \mset{K}^{\times}$ is estimated as the proportion of the data points such that the uncertainty part related to $a$ is attributed to cluster $k_{a}$ for all $a \in \mset{A}$.
To take into account the weak correlation between groups of attributes (so far assumed to be independent),
we retain only the most probable combinations of group-scenarios
$k \mdef \pab{k_{a} \in \mset{K}_{a}, a \in \mset{A}} \in \mset{K}^{\times}$,
keeping those with associated probability greater than a threshold probability~$\tilde{p}$.
The more conservative we desire to be,
the greater the acceptance probability threshold~$\tilde{p}$ should be.
The set of scenarios~$\mset{K}$ is then constituted of the combinations in~$\mset{K}^{\times}$ satisfying the probability threshold~$\tilde{p}$.
The probability of occurrence~$p_{k}$ of each scenario~$k \in \mset{K}$ used in~\eqref{eq:RO_objective} is then defined as the scaled probability when the scenarios that do not satisfy the acceptance threshold are excluded,
\ie $p_{k} = \hat{p}_{k} / \sum_{k \in \mset{K}} \hat{p}_{k}$.

Finally,
the PUS in the original basis related to a scenario $k = \pab{k_{a} \in \mset{K}_{a}, a \in \mset{A}} \in \mset{K}$ is then defined as follows:
\begin{equation}
  \label{eq:polyhedral_uncertainty_set_scenario}
  \mset{V}_{k}
  \mdef
  \setc{
    \vb \in \R^{24n}
    }{
      \forall a \in \mset{A}:
      \begin{cases}
      \bar{\wb}_{a} \in \bar{\mset{W}}_{a,k_{a}}^{pol} \\
      \vb_{a} = \Pb_{a}^{T}
      \begin{bmatrix}
	\bar{\wb}_{a} \\  \0_{n_{a} - r_{a}}
      \end{bmatrix}
      \end{cases}
  }.
\end{equation}
where $\0_{n_{a} - r_{a}}$ is a vector of zeros of size $n_{a} - r_{a}$.
The reduction in dimension appears clearly:
the uncertain part $\vb_{a} \in \R^{24n_{a}}$ related to~$a \in \mset{A}$ is generated from elements of $\bar{\mset{W}}_{a,k_{a}}^{pol}$,
which is in the lower-dimensional space~$\R^{r_{a}}$.

\paragraph{Uncertain set for several successive days.}%
\label{par:successive_days}

In some applications (\eg unit commitment),
we are interested in taking into account the uncertainty over several days.
Given the succession of~$N$ daily scenarios $k_{1 \to N} \mdef \pa{k_{1}, \dots, k_{N}} \in \mset{K} \times \dots \times \mset{K}$,
we can define an associated PUS based on the definition of $\mset{V}_{k}$ in~\eqref{eq:polyhedral_uncertainty_set_scenario} as follows:
{\small
\begin{equation}
  \label{eq:polyhedral_uncertainty_set_scenario_successive_days}
  \mset{V}_{k_{1 \to N}}
  \mdef
  \setc{
    \vb \in \R^{24nN}
    }{
      \begin{cases}
    \vb = \pac{\vb_{1}, \dots, \vb_{N}} \\
	\forall i \in \rangeo{N}:
	\vb_{i} \in \mset{V}_{k_{i}}\\
	\forall i \in \rangeo{N-1}:
      \abs{ \spa{\Vb_{i+1}[:,1] - \Vb_{i}[:,24]} - \mub_{\Delta}} \circ \frac{1}{\sigmab_{\Delta}} \leq c
      \end{cases}
  },
\end{equation}
}%
where $\frac{1}{\sigmab_{\Delta}}$ must be understood componentwise,
where $\Vb_{i} \in \R^{n \times 24}$ is the matrix representation of~$\vb_{i}$,
and where
$\mub_{\Delta}$ and $\sigmab_{\Delta}$ are vectors in $\R^{n}$ that correspond respectively to the sample mean and standard deviation of the difference between the first value and last value of the previous day in the daily profile for each~$n$ attributes.  
The third constraint in~\eqref{eq:polyhedral_uncertainty_set_scenario_successive_days}
enforces continuity between the daily profiles with the parameter $c > 0$ as it limits the step size between successive profiles.

\newpage
\section{Case study: uncertainty related to energy data in the UK}%
\label{sec:numerical_example}

\begin{wrapfigure}{r}{0.6\textwidth}
  \includegraphics[width=1.0\linewidth]{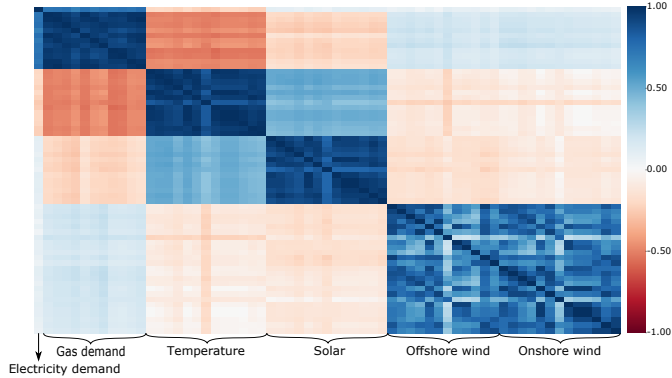}
  \caption{Correlation between all the attributes}%
  \label{fig:correlation_stage0}
  \vspace{-0.4cm}
\end{wrapfigure}
We apply the method presented in the previous section to derive uncertainty sets for quantifiable sources of uncertainty related to power generation based on historical data for the UK in 2015.
This includes national electricity demand, as well as regional data on gas demand, temperature, solar availability, offshore and onshore wind availability
(we refer to the paper \citet{Charitopoulos2022a} for more details).
Each data point corresponds to the collection of hourly daily profiles,
of each of the previously mentioned attributes.
The size of the uncertainty space in which each data point resides is therefore high dimensional,
which motivates the need for a reduction technique to generate uncertainty sets.

Based on the correlation plot
(see \cref{fig:correlation_stage0}),
we observe three groups of attributes that are strongly correlated:
\emph{seasonal} data (electricity/gas demand and temperature),
\emph{solar} data
and
\emph{wind} data (off/onshore wind),
\ie $\mset{A} = \pab{seasonal, solar, wind}$.

\begin{wrapfigure}{r}{0.4\textwidth}
  \includegraphics[width=1.0\linewidth]{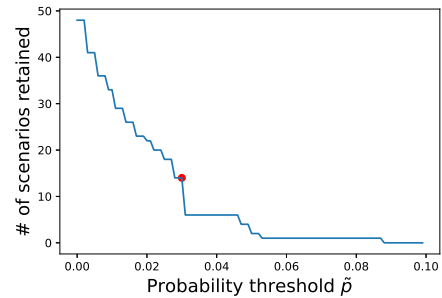}
  \caption{
      Number of scenarios given the probability threshold $\tilde{p}$.
  }%
  \label{fig:p_number_scenarios}
  \vspace{-0.4cm}
\end{wrapfigure}

\begin{figure}[h!]
  \centering
  \hfill
  \subfloat{{\includegraphics[width=0.33\linewidth]{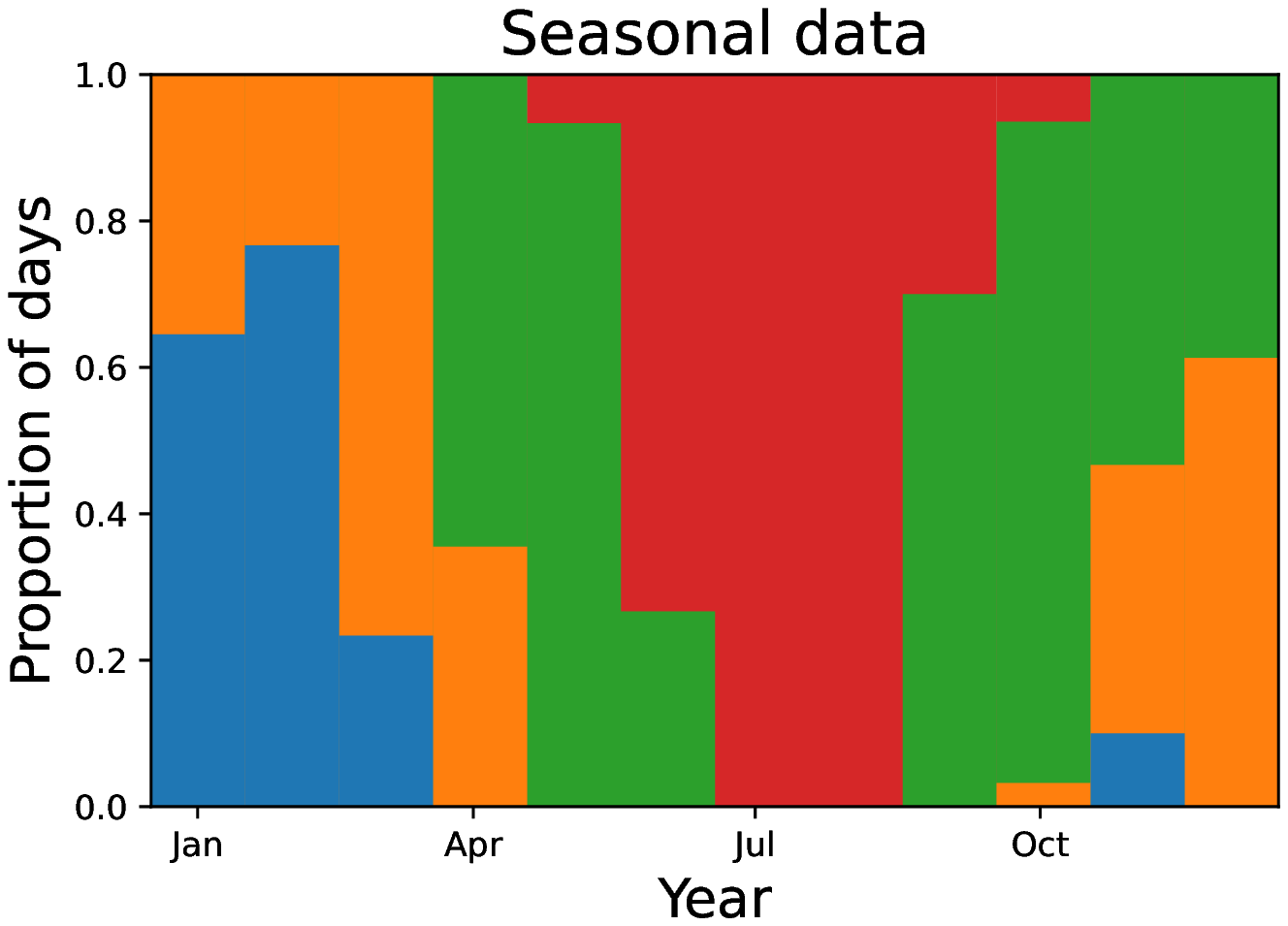}}}
  \hfill
  \subfloat{{\includegraphics[width=0.33\linewidth]{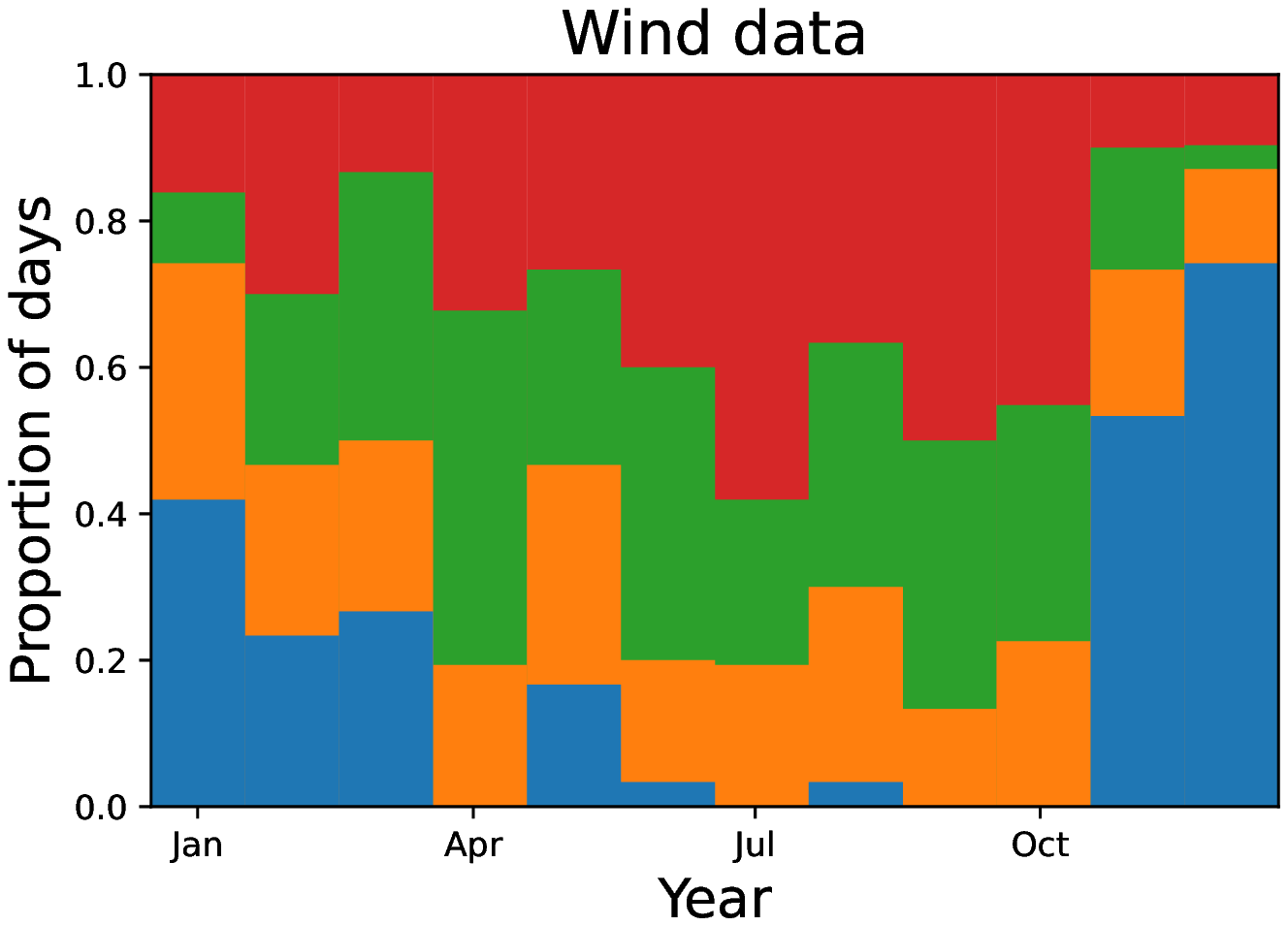}}}
  \hfill
  \subfloat{{\includegraphics[width=0.33\linewidth]{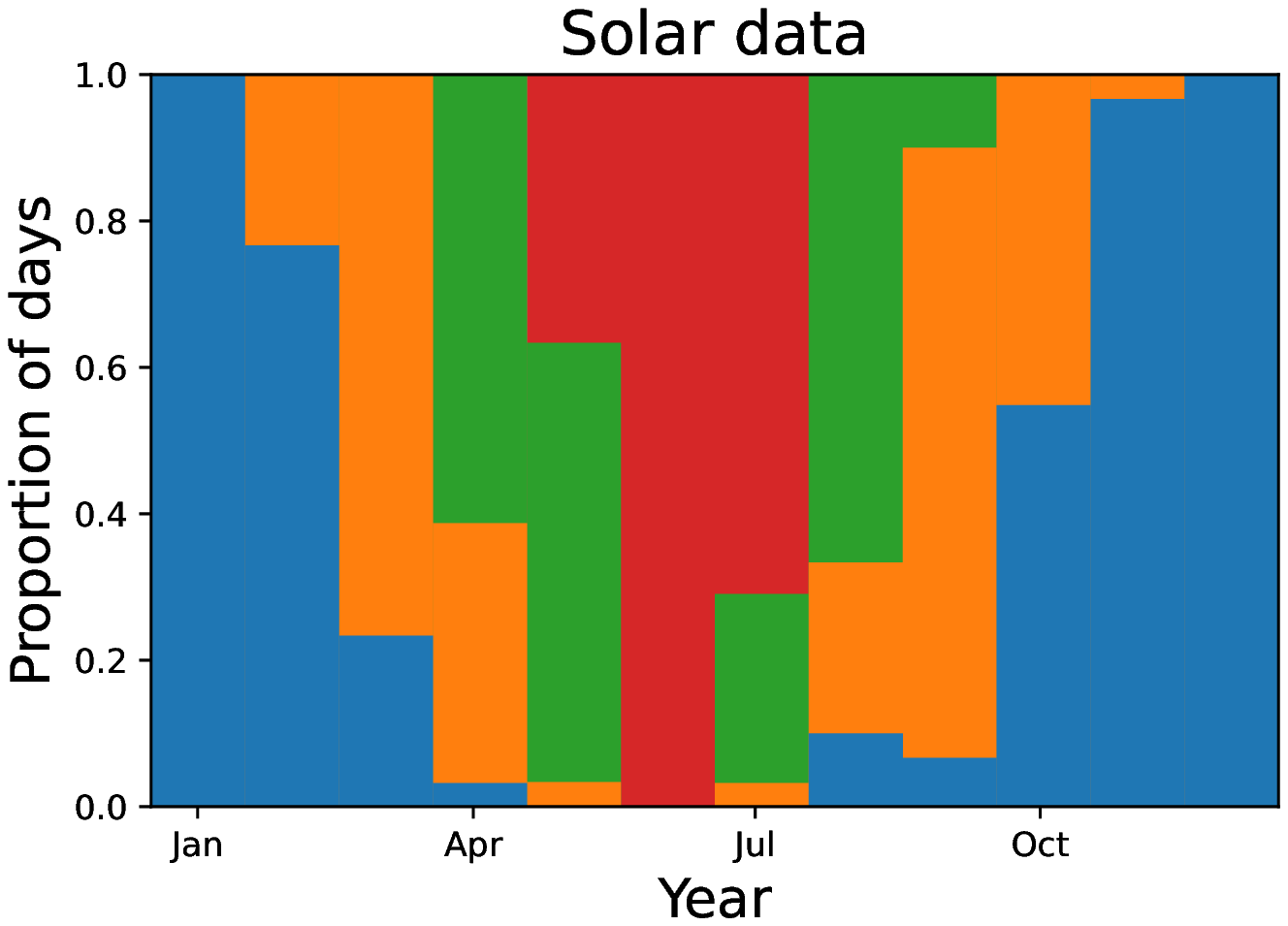}}}
  \hfill
  \subfloat{{}}
  \hfill
  \caption{
      Proportion of days allocated to each cluster for all groups of attributes
  }%
  \label{fig:clustering_histogram}
\end{figure}

\Cref{fig:clustering_histogram} illustrates for each month the proportions of days that are allocated to each cluster (4 clusters are generated for each group of attributes $a \in \mset{A}$).
The correlation between the solar and seasonality data is clearly apparent:
the simultaneous occurrence of the blue wind cluster with the red solar cluster is for instance highly improbable. 
Now,
if we focus on the clustering from April to October,
we observe that the clusters proportions in the wind data are roughly the same for each month.
This translates the idea that the seasonal data and wind data are independent for this period of the year.
As a consequence of these observations,
we do not define the scenarios as the Cartesian product of group-scenarios since all combinations are not as probable.
\Cref{fig:p_number_scenarios} depicts how the number of scenarios retained varies with the acceptance threshold~$\tilde{p}$.
It seems that adopting a probability threshold of $\tilde{p} = 3\%$,
which results in keeping 14 scenarios,
is in line with the trade-off between accuracy and computational tractability:
a higher value would lead to too few scenarios (6 or less),
while a lower value would lead to too many scenarios and therefore would not allow for the pruning of the least likely realisations.

\begin{figure}[h!]
  \centering
  \hfill
  \subfloat{{\includegraphics[width=0.45\linewidth]{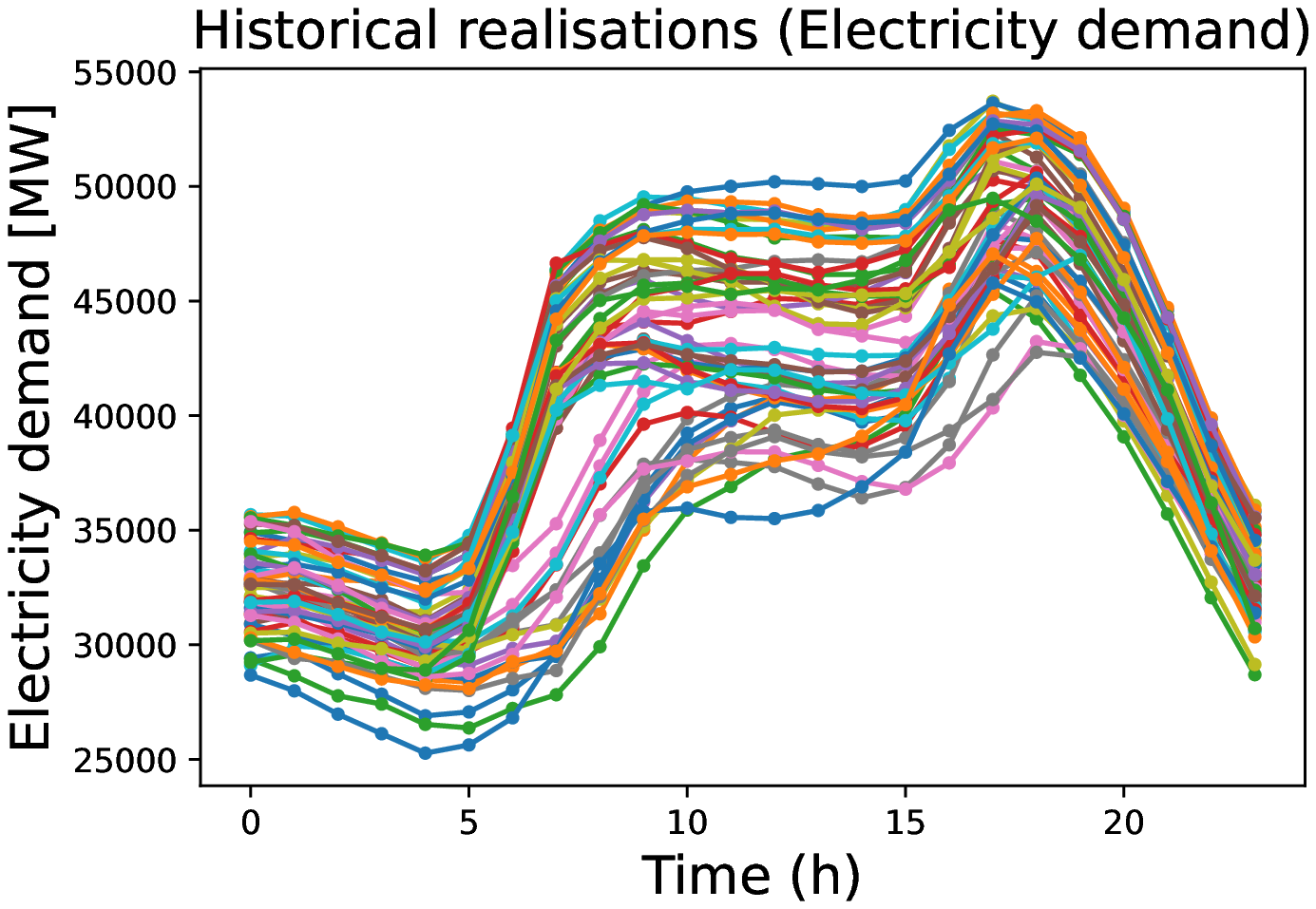}}}
  \hfill
  \subfloat{{\includegraphics[width=0.45\linewidth]{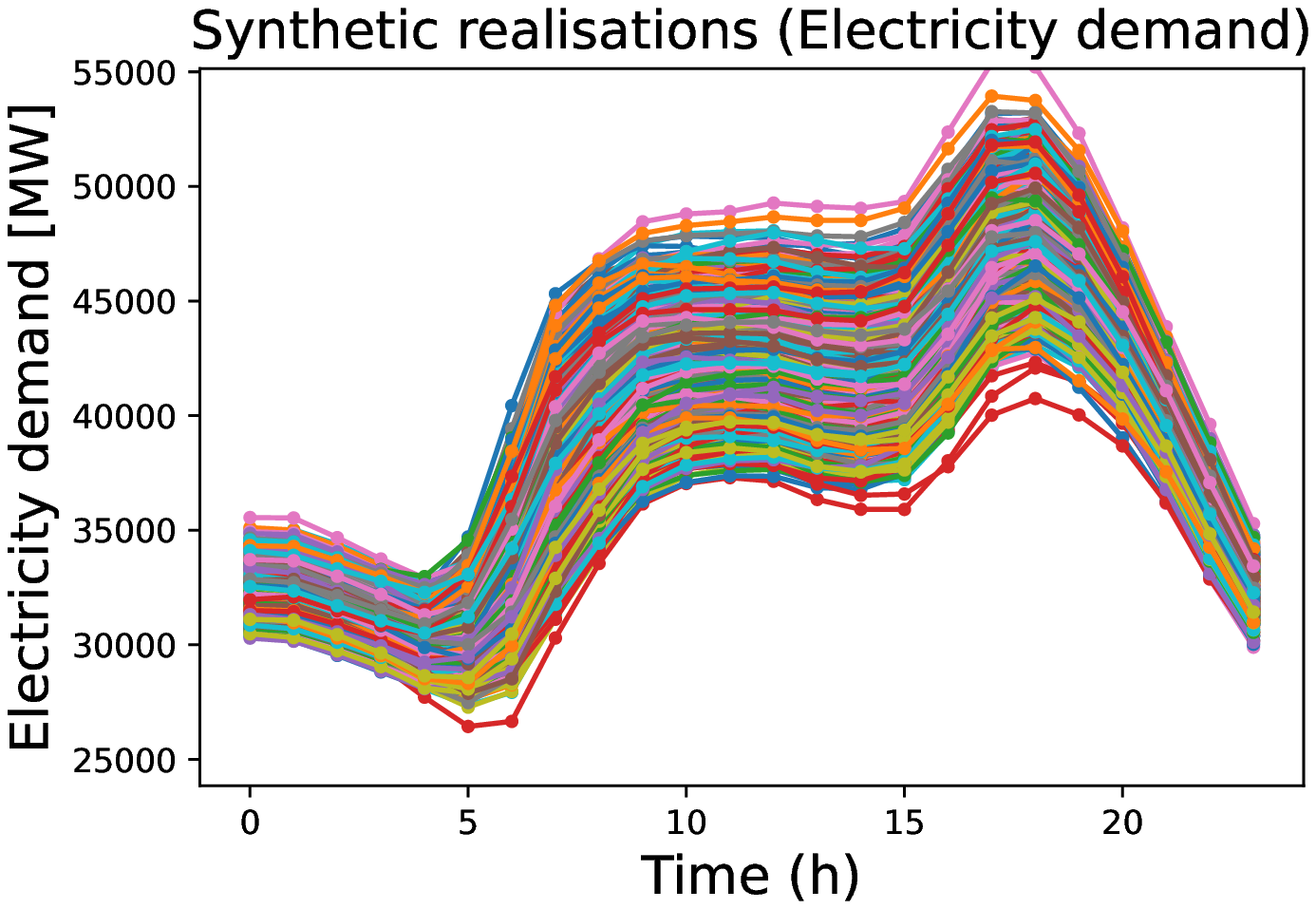}}}
  \hfill 
  \\
  \hfill
  \subfloat{{\includegraphics[width=0.45\linewidth]{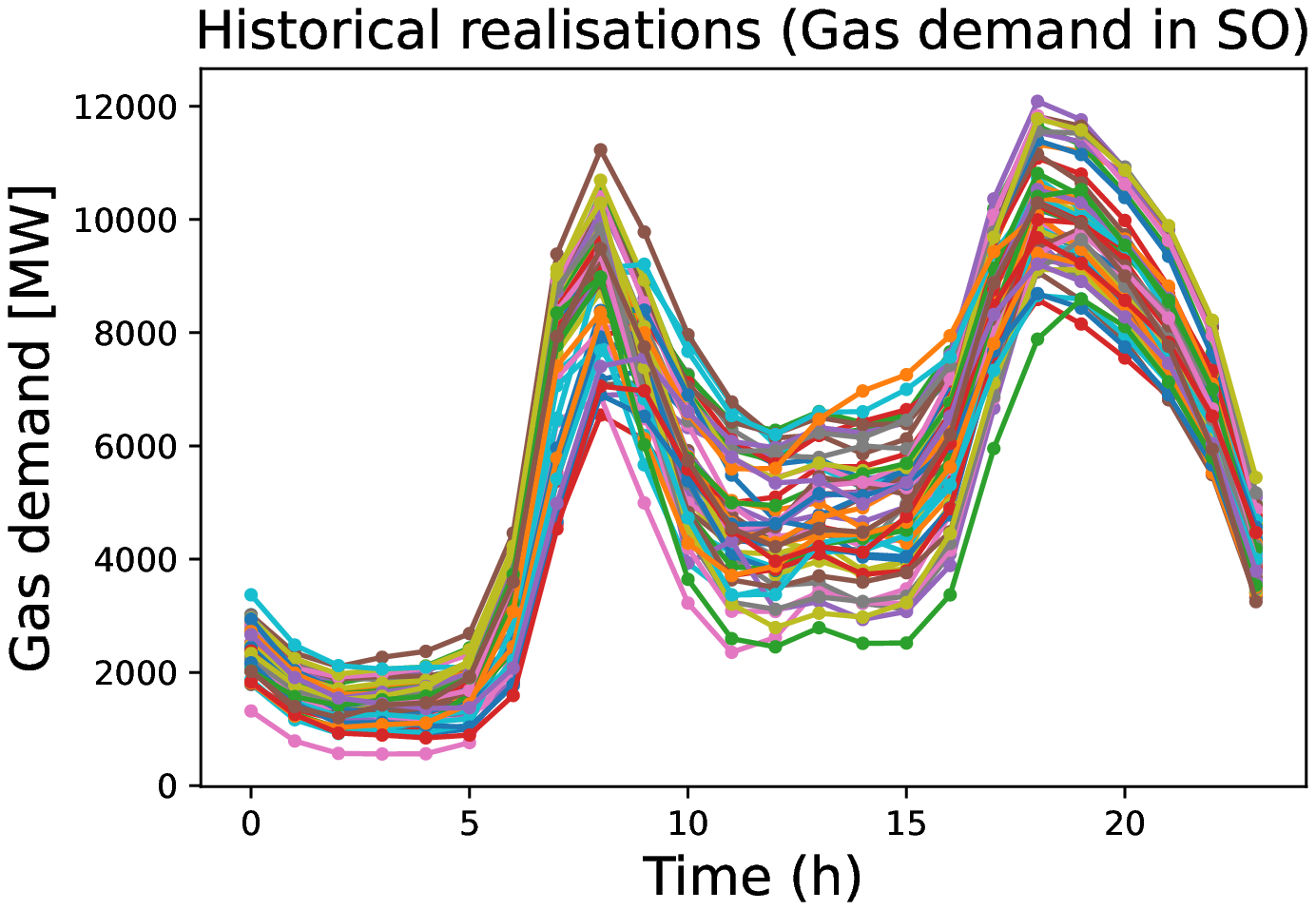}}}
  \hfill
  \subfloat{{\includegraphics[width=0.45\linewidth]{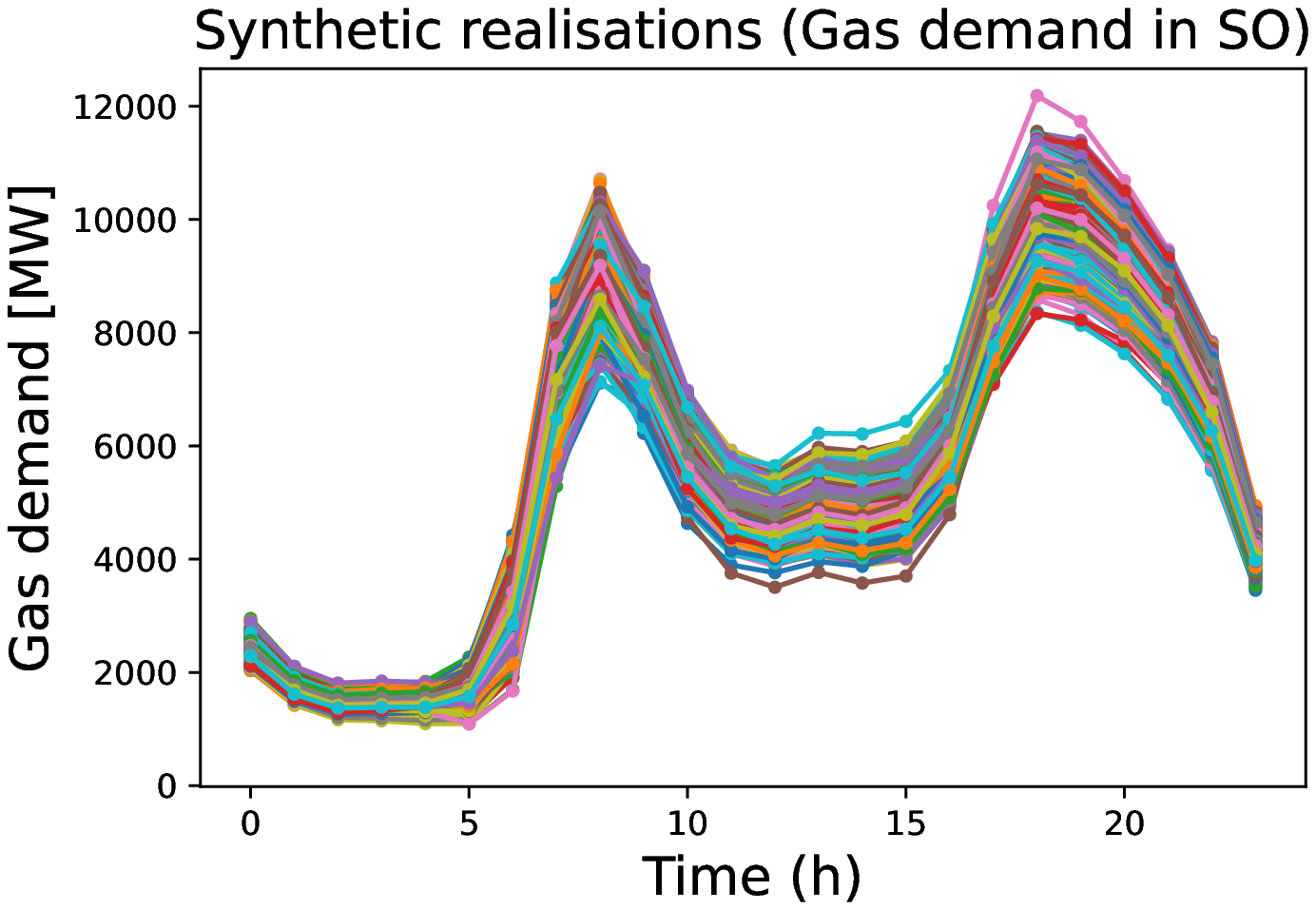}}}
  \hfill
  \caption{
  Comparison of historical and synthetic realisations for the blue cluster.
}
  \label{fig:pus_cluster}
\end{figure}

\begin{figure}[h!]
  \centering
  \includegraphics[width=1.0\linewidth]{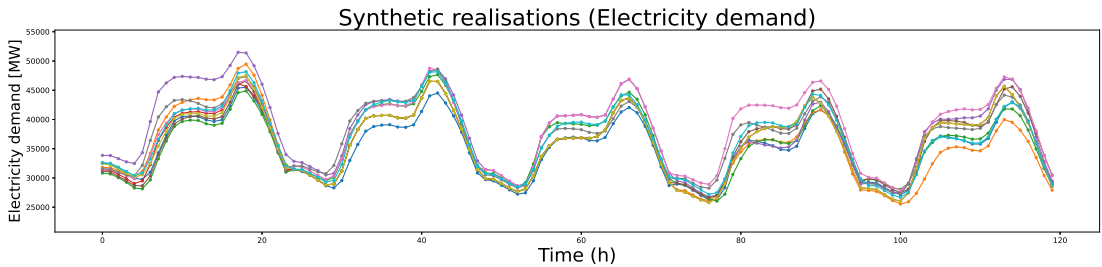}
  \caption{
Synthetic realisations for 
$k_{1 \to 5} = \pac{blue, blue, orange, orange, orange}$.
}
  \label{fig:successive_clusters}
\end{figure}

Finally,
\cref{fig:pus_cluster} depicts data generated with~\eqref{eq:polyhedral_uncertainty_set_scenario}
for the parameters values:
$a = seasonal$, $k_{a} = blue$, $r_{a} = 24$, $\alpha_{a,k_{a}} = 5\%$, $\Phi_{a,k_{a}} = 5$, $s_{a,k_{a}} = 24 $ and $\Psi_{a,k_{a}} = 1.5$.
We observe that the daily outliers in the historical data are ignored and,
as desired,
that the synthetic realisations are similar to the 90\% of cases around the median. 
The size reduction is as follows: $n_{a} = 40$ attributes are associated with the group $a = seasonal$,
which implies that the seasonal data is of size $24n_{a} = 960$ and is generated on the basis of $r_{a} = 24$ PCA components;
this corresponds to a size reduction factor of 40.
\Cref{fig:successive_clusters} depicts realisations generated with~\eqref{eq:polyhedral_uncertainty_set_scenario_successive_days} for the seasonal group of attributes for the succession of five days,
where $k_{1 \to 5} = \pac{blue, blue, orange, orange, orange}$
and where the continuity constraint parameter is equal to $c = 2.5$.


\section{Results and discussion}%
\label{sec:results}

We first point out why we did not start by clustering the data and then, for each cluster,
generate a PUS,
as it would implicitly take into account the (weak) correlation between attributes and avoid generating clusters based on a correlation analysis.
The answer lies in two main points: 
Firstly, the scarcity of data.
Indeed, in the case study,
if we were to generate the 14 scenarios by clustering over the 365 historical data points,
then each of them would have on average 26 data points.
In contrast,
in the method presented in this paper,
only 4 clusters are generated per group of attributes. 
The quantiles of each cluster are then derived on the basis of around 90 data points.
The PUS derived with our method are thus less subject to statistical estimate errors. 
Secondly,
establishing groups of uncorrelated attributes allows for associating to each one of them a uncertainty dispersion budget~$\Phi_{a,k_{a}}$. 
In this way,
we do not constrain the realisations of independent attributes together.
\section{Concluding remarks}%
\label{sec:conclusion}
In conclusion, the advantages of our method are the following:
First, it considerably reduces the size of the original uncertainty set as it can be expressed as a linear transformation of a truncated basis using PCA.
Secondly,
the PUS generated are statistically meaningful and can exploit the possibly independence of uncertain attributes.
Thirdly, it generates consistent uncertainty daily profiles for which the continuity from hour to hour is guaranteed.
Finally, it derives polyhedral uncertainty sets that can then be used in a robust optimisation (RO) type of problem.

\textbf{Reproducibility:}
All the graphs of this article can be reproduced with the code made available on \href{https://github.com/julien-vaes/A_data_driven_uncertainty_modelling_approach_for_energy_optimisation_problems.git}{\texttt{GitHub}},
where we have also shared the 2015 UK energy uncertainty dataset.

\textbf{Acknowledgements:}
Financial support from the EPSRC (under project EP/T022930/1) is gratefully acknowledged.



{\small\bibsep=0pt
\bibliography{latex_files/references.bib}

\begin{thebibliography}{12}
\providecommand{\natexlab}[1]{#1}
\providecommand{\url}[1]{\texttt{#1}}
\expandafter\ifx\csname urlstyle\endcsname\relax
  \providecommand{\doi}[1]{doi: #1}\else
  \providecommand{\doi}{doi: \begingroup \urlstyle{rm}\Url}\fi

\bibitem[Ben-Tal et~al.(2009)Ben-Tal, {El Ghaoui}, and Nemirovski]{Ben-Tal2009}
A.~Ben-Tal, L.~{El Ghaoui}, and A.~Nemirovski.
\newblock \emph{{Robust Optimization}}.
\newblock Princeton University Press, 2009.
\newblock \doi{10.1515/9781400831050}.

\bibitem[Bertsimas and Sim(2004)]{Bertsimas2004a}
D.~Bertsimas and M.~Sim.
\newblock {The Price of Robustness}.
\newblock \emph{Oper. Res.}, 52\penalty0 (1):\penalty0 35--53, 2004.
\newblock \doi{10.1287/opre.1030.0065}.

\bibitem[Bianchi et~al.(2018)Bianchi, Gupta, and Kallus]{Bertsimas2018}
M.~Bianchi, V.~Gupta, and N.~Kallus.
\newblock {Data-driven robust optimization}.
\newblock \emph{Math. Program.}, 167\penalty0 (2):\penalty0 235--292, 2018.
\newblock \doi{10.1007/s10107-017-1125-8}.

\bibitem[Birge and Louveaux(2011)]{Birge2011}
J.~R. Birge and F.~Louveaux.
\newblock \emph{{Introduction to Stochastic Programming}}.
\newblock Springer New York, 2011.
\newblock \doi{10.1007/978-1-4614-0237-4}.

\bibitem[Charitopoulos et~al.(2022)Charitopoulos, Fajardy, Chyong, and
  Reiner]{Charitopoulos2022a}
V.~M. Charitopoulos, M.~Fajardy, C.~K. Chyong, and D.~Reiner.
\newblock {The case of 100{\%} electrification of domestic heat in Great
  Britain}.
\newblock 2022.
\newblock \doi{10.17863/CAM.81913}.

\bibitem[Chen(2017)]{Chen2017}
Y.~C. Chen.
\newblock {A tutorial on kernel density estimation and recent advances}.
\newblock \emph{Biostat. Epidemiol.}, 1\penalty0 (1):\penalty0 161--187, 2017.
\newblock \doi{10.1080/24709360.2017.1396742}.

\bibitem[Ding and He(2004)]{Ding2004}
C.~Ding and X.~He.
\newblock {K -means clustering via principal component analysis}.
\newblock In \emph{Twenty-first Int. Conf. Mach. Learn. - ICML '04}, page~29.
  ACM Press, 2004.
\newblock \doi{10.1145/1015330.1015408}.

\bibitem[Ning and You(2017)]{Ning2017}
C.~Ning and F.~You.
\newblock {Data‐driven adaptive nested robust optimization: General modeling
  framework and efficient computational algorithm for decision making under
  uncertainty}.
\newblock \emph{AIChE J.}, 63\penalty0 (9):\penalty0 3790--3817, 2017.
\newblock \doi{10.1002/aic.15717}.

\bibitem[Ning and You(2018)]{Ning2018}
C.~Ning and F.~You.
\newblock {Data-driven decision making under uncertainty integrating robust
  optimization with principal component analysis and kernel smoothing methods}.
\newblock \emph{Comput. Chem. Eng.}, 112:\penalty0 190--210, 2018.
\newblock \doi{10.1016/j.compchemeng.2018.02.007}.

\bibitem[Ning and You(2019)]{Ning2019a}
C.~Ning and F.~You.
\newblock {Data-Driven Adaptive Robust Unit Commitment Under Wind Power
  Uncertainty: A Bayesian Nonparametric Approach}.
\newblock \emph{IEEE Trans. Power Syst.}, 34\penalty0 (3):\penalty0 2409--2418,
  2019.
\newblock \doi{10.1109/TPWRS.2019.2891057}.

\bibitem[Roald et~al.(2023)Roald, Pozo, Papavasiliou, Molzahn, Kazempour, and
  Conejo]{Roald2023a}
L.~A. Roald, D.~Pozo, A.~Papavasiliou, D.~K. Molzahn, J.~Kazempour, and
  A.~Conejo.
\newblock {Power systems optimization under uncertainty: A review of methods
  and applications}.
\newblock \emph{Electr. Power Syst. Res.}, 214:\penalty0 108725, 2023.
\newblock \doi{10.1016/j.epsr.2022.108725}.

\bibitem[Shang et~al.(2017)Shang, Huang, and You]{Shang2017}
C.~Shang, X.~Huang, and F.~You.
\newblock {Data-driven robust optimization based on kernel learning}.
\newblock \emph{Comput. Chem. Eng.}, 106:\penalty0 464--479, 2017.
\newblock \doi{10.1016/j.compchemeng.2017.07.004}.

\end{thebibliography}
}

\end{document}